\documentclass[12pt]{amsart}
\usepackage{amsmath}
\usepackage{amssymb}
\usepackage{amsthm}
\usepackage{geometry}
\usepackage{hyperref}
\usepackage{enumitem}
\usepackage{amsfonts}
\usepackage{amscd}
\usepackage{relsize}
\usepackage{setspace}
\usepackage{url}
\usepackage{xspace}
\usepackage{mathrsfs}
\usepackage{dsfont}
\usepackage{lmodern}

\usepackage[textsize=footnotesize,textwidth=15ex]{todonotes}

\usepackage{xcolor}

\definecolor{forestgreen}{rgb}{0.13, 0.55, 0.13}

\newtheorem{theorem}{Theorem}[section]
\newtheorem{lemma}[theorem]{Lemma}
\newtheorem{prop}[theorem]{Proposition}
\newtheorem{corollary}[theorem]{Corollary}

\theoremstyle{definition}
\newtheorem{remark}[theorem]{Remark}

\newtheorem{definition}[theorem]{Definition}
\newtheorem{question}[theorem]{Question}

\newcommand{\triv}{\{1\}}

\newcommand{\QZ}{\mathrm{QZ}}

\newcommand{\CC}{\mathrm{C}}

\newcommand{\N}{\mathrm{N}}

\newcommand{\tdlc}{t.d.l.c.\@\xspace}

\newcommand{\ldlat}{\mathcal{LD}}

\newcommand{\Ends}{\mathscr{E}}
\newcommand{\rist}{\mathrm{rist}}

\newcommand{\lnorm}{\mathcal{LN}}

\newcommand{\lcent}{\mathcal{LC}}

\newcommand{\defbold}{\textbf}

\newcommand{\bN}{\mathbb{N}}

\newcommand{\bQ}{\mathbb{Q}}
\newcommand{\bR}{\mathbb{R}}
\newcommand{\bZ}{\mathbb{Z}}

\newcommand{\mc}{\mathcal}
\newcommand{\mf}{\mathfrak}
\newcommand{\ms}{\mathscr}

\newcommand{\inv}{^{-1}}

\DeclareMathOperator{\Aut}{Aut}

\newcommand{\cgrp}[1]{\overline{\langle #1 \rangle}}

\newcommand{\grp}[1]{\langle #1 \rangle}
\newcommand{\ol}[1]{\overline{#1}}

\begin{document}

\title{Growing trees from compact subgroups}
\thanks{TM: F.R.S.-FNRS Research associate; PEC and TM are supported in part by the  F.R.S.-FNRS and the FWO under the EOS programme (project ID 40007542); CDR: Research supported in part by ARC grant FL170100032.}

\author[P.-E. Caprace]{Pierre-Emmanuel \textsc{Caprace}}
\address{Pierre-Emmanuel Caprace\\
UCLouvain, IRMP\\
Chemin du Cyclotron 2, Box L7.01.02, 1348 Louvain-la-Neuve, Belgium}
\email{pierre-emmanuel.caprace@uclouvain.be}

\author[T. Marquis]{Timoth\'ee \textsc{Marquis}}
\address{Timoth\'ee Marquis\\
UCLouvain, IRMP\\
Chemin du Cyclotron 2,  Box L7.01.02, 1348 Louvain-la-Neuve, Belgium}
\email{timothee.marquis@uclouvain.be}

\author[C. Reid]{Colin D. \textsc{Reid}}
\address{Colin D. Reid\\
The University of Newcastle, School of Information and Physical Sciences\\
Callaghan, NSW 2308, Australia}
\email{colin@reidit.net}

\begin{abstract}
We establish a new connection between local and large-scale structure in compactly generated totally disconnected locally compact (t.d.l.c.) groups $G$, finding a sufficient condition for $G$ to have more than one end in terms of its compact subgroups.  The condition actually results in an action of a quotient group $G/N$ on a tree with faithful micro-supported action on the boundary, where $N$ is compact, and is closely related to the Boolean algebra formed by the centralisers of the subgroups of $G/N$ with open normaliser.  As an application, we find a sufficient condition, given a one-ended t.d.l.c. group $G$, for all direct factors of open subgroups of $G$ to be trivial or open.
\end{abstract}

\maketitle

\setcounter{tocdepth}{1}

\section{Introduction}

\begin{flushright}
	\begin{minipage}[t]{0.55\linewidth}\itshape\small
This is not our world with trees in it. It's a world of trees, where humans have just arrived.
		
		\hfill\upshape (Richard Powers, \emph{The Overstory},  2018)
	\end{minipage}
\end{flushright}

A recurring theme in the theory of totally disconnected, locally compact (\tdlc) groups is the relationship between local structure, that is, properties evident in an arbitrarily small neighbourhood of the identity, and global or large-scale properties, for instance in the sense of quasi-isometry invariants of the group, or ``commability'' in the sense of \cite{CorCommable}.  The connection between local structure and large-scale properties is looser than for Lie groups, but there are still non-trivial relationships between the two.  The goal of this article is to investigate a certain aspect of this relationship for general compactly generated \tdlc groups: we provide an algebraic criterion on the local structure of a \tdlc group ensuring that the group has infinitely many ends and, hence, that it acts on a tree with compact open edge stabilizers.  An application is given in another article by the same authors, showing that many one-ended locally compact Kac--Moody groups are locally indecomposable: see \cite{LocindecKM}.

A fundamental example of local structure is the \defbold{structure lattice} $\lnorm(G)$ of a \tdlc group, introduced in \cite{CRWpart1}.  This is the poset $\lnorm(G) = \mathrm{LN}(G)/\sim_o$, where $\mathrm{LN}(G)$ is the set of closed \defbold{locally normal} subgroups (subgroups with open normaliser) of $G$ ordered by inclusion and ${H \sim_o K}$ if $H \cap K$ is open in $H$ and $K$. Note that $\lnorm(G)$ is equipped with an action of $G$ by conjugation, which we also carry over to conjugation-invariant subsets of $\lnorm(G)$.  To avoid some complications it is useful to assume that $G$ is \defbold{[A]-semisimple}, meaning that the subgroup $\QZ(G)$ of elements with open centraliser is trivial and that $G$ has no nontrivial abelian locally normal subgroups.  We set
\[
\mathrm{LC}(G) := \{ \CC_G(K) \mid K \in \mathrm{LN}(G)\}; \; \lcent(G) := \mathrm{LC}(G)/\sim_o.
\]
When $G$ is [A]-semisimple, then $\mathrm{LC}(G)$ is a Boolean algebra called the \defbold{(global) centraliser lattice} of $G$; its elements are pairwise inequivalent under $\sim_o$, so it is isomorphic as a Boolean algebra to the \defbold{(local) centraliser lattice} $\lcent(G)$.  There is also a subalgebra $\ldlat(G)$ of $\lcent(G)$, the \defbold{local decomposition lattice}, which consists of those elements of $\lnorm(G)$ represented by direct factors of open subgroups.  In particular, we say $G$ is \defbold{locally indecomposable} if $\ldlat(G)$ is trivial and \defbold{faithful locally decomposable} if $G$ acts faithfully on $\ldlat(G)$.  When $G$ is [A]-semisimple, having faithful action on $\lcent(G)$ is equivalent to the existence of a faithful action on a compact totally disconnected Hausdorff space $X$ that is \defbold{micro-supported}, meaning that for every nondense subset $Y$ of $X$, there is $g \in G \setminus \{1\}$ that fixes $Y$ pointwise.

An interesting class to consider in this context is the class $\ms{S}$ of nondiscrete, compactly generated, topologically simple \tdlc groups.  Unlike with simple Lie groups, the local structure does not determine the global structure: S. Smith (\cite{SmithDuke}) exhibited a family of $2^{\aleph_0}$ pairwise nonisomorphic groups in $\ms{S}$ that are locally isomorphic to one another and all faithful locally decomposable.  However, it is the case that all groups $G \in \ms{S}$ are [A]-semisimple, and that $G$ acts faithfully on any $G$-invariant subalgebra of $\lcent(G)$ other than the trivial one $\{0, \infty\}$ (there could be none).

The space $\Ends G$ of \textbf{ends} of a locally compact group $G$ is a large-scale invariant; see \cite[Definition~8.B.12]{CdH16} \S\ref{subsec:ends} for the precise definition of $\Ends G$ in use in this article. A result of H.~Abels \cite{Abels-74} ensures that an analogue of Stallings' splitting theorem holds for compactly generated \tdlc groups: such a group has infinitely many ends if and only if it has a continuous  unbounded action on a tree with compact open edge stabilizers (see also Lemma~\ref{lem:tree_geometric} and Proposition~\ref{prop:end_types} below). Moreover, $2$-ended groups have a restricted and well understood algebraic structure, from which it follows that a group in $G \in \ms{S}$ has either one end or infinitely many ends. 

To avoid ambiguity, we will refer to edges of graphs as \defbold{arcs} when taking account of orientation, \textit{i.e.} the arc from $v$ to $w$ is distinct from the arc from $w$ to $v$.  Many examples of groups in $\ms{S}$, including Smith's examples, have been constructed as groups acting \defbold{arc-geometrically} on a leafless tree, meaning with compact open arc stabilisers and preserving no proper subtree.  As mentioned above, for a compactly generated \tdlc group, having such an action is equivalent to having more than one end.  Write $\ms{S}_{\infty}$ for the infinitely-ended groups in $\ms{S}$ and $\ms{S}_{\mc{LD}}$ for the locally decomposable groups in $\ms{S}$.  Neither of $\ms{S}_{\infty}$ and $\ms{S}_{\mc{LD}}$ contains the other: for example, $\mathrm{PSL}_2(\bQ_p)$ is infinitely-ended but locally indecomposable, whereas Neretin's groups are one-ended (see \cite[Corollary~9.0.12]{Kerg}; an alternative proof of that fact can be deduced from the main result of \cite{Genevois}) but locally decomposable.  There are also groups in $\ms{S}$ that are in neither $\ms{S}_{\infty}$ nor $\ms{S}_{\mc{LD}}$.  For example, given a simple algebraic group $G$ over a local field $k$, the associated group $S = G(k)^+/Z \in \ms{S}$ is locally indecomposable (see \cite[Appendix~A]{CRWpart2}).  However $S$ acts geometrically on a building $X$ with apartments isometric to $\bR^\ell$, where $\ell$ is the $k$-rank of $G$ (\cite{BruhatTits}); if $\ell \ge 2$ then the apartments are one-ended, so $X$ is one-ended (\textit{e.g.} \cite[Lemma~5.21]{LocindecKM}) and hence $S$ is one-ended.

Nevertheless, there is a large overlap between $\ms{S}_{\infty}$ and $\ms{S}_{\mc{LD}}$ in the known examples, and for certain geometric constructions they have been shown to be equivalent.  An example of such a construction occurs in \cite{CapraceRAB}, where the equivalence between having infinitely many ends and being locally decomposable is non-vacuous, in the sense that some of the groups in $\ms{S}$ obtained in \cite{CapraceRAB} are one-ended and not locally decomposable, while others are infinitely-ended and locally decomposable.  More generally, except in the setting of algebraic groups over local fields, it is difficult to show a group $G$ in $\ms{S}$ is locally indecomposable, and even more difficult to show $G$ has no micro-supported action on the Cantor space.  In fact, the majority of constructions to date of groups in $\ms{S}$ are based on J. Tits' property (P) (\cite{Tits70}) or its generalisations based on the double commutator lemma (see \cite[Proposition~I]{CRWpart2}), which often directly imply local decomposability, or at least a micro-supported action.

In the present article, we establish a sufficient condition for a compactly generated \tdlc group $G$ to act arc-geometrically on a tree $T$, such that the action on an associated compact boundary $\ol{\Ends} T$ (see Section~\ref{subsec:ends}) is micro-supported.  To simplify the presentation here we assume that the only compact normal subgroup of $G$ is the trivial subgroup $\triv$; for the full statement without this assumption, see Section~\ref{subsec:cpctend}.

\begin{definition}
Let $G$ be a \tdlc group without nontrivial compact normal subgroup and let $K$ be an infinite compact subgroup of $G$.  We say $K$ is a \defbold{TMS subgroup} of $G$ if there is a compact open subgroup $U$ of $G$ with the following properties:
\begin{enumerate}[label=(\alph*)]
\item For all conjugates $V$ of $U$ in $G$, the set of $g \in G$ such that $gKg\inv \le U$ but $gKg\inv \nleq V$ has compact closure.
\item The set of $g \in G$ such that $gKg\inv \le U$ does not have compact closure.
\end{enumerate}
\end{definition}

The letters TMS stand for \textit{tree micro-supported}.  The definition and choice of terminology are motivated by the following results, where a \defbold{half-tree} of a tree $T$ is a proper subtree of $T$ that is joined to the rest of $T$ by a single edge.

\begin{prop}[See Proposition~\ref{prop:smooth_tree}]\label{intro:tree->TMS}
Let $G$ be a compactly generated \tdlc group acting faithfully and arc-geometrically on a leafless tree $T$.  Suppose that for some half-tree $T'$ of $T$, the pointwise fixator $K$ of $T'$ in $G$ fixes only finitely many arcs of $T \setminus T'$.  Then $K$ is a TMS subgroup of $G$.
\end{prop}

\begin{theorem}[See Theorem~\ref{cpctend}]\label{intro:cpctend}
Let $G$ be a compactly generated \tdlc group without nontrivial compact normal subgroup.  Suppose that $G$ has a TMS subgroup $K$.  Then $G$ is [A]-semisimple.  Moreover, $G$ acts faithfully on a leafless tree $T$, such that the pointwise fixator of every half-tree contains a conjugate of $K$ and such that one of the following holds:
\begin{enumerate}[label=(\roman*)]
\item $T$ is locally finite and $G$ acts vertex-transitively with compact open stabilisers on $T$, fixing a unique end;
\item $G$ preserves no end or proper subtree of $T$ and $G$ has compact open arc stabilisers on $T$.
\end{enumerate}
In either case, the action of $G$ on $\Ends G$ is faithful and nondiscretely micro-supported.
\end{theorem}

We remark that in case (i), the group $G$ is a \textit{focal hyperbolic group} in the sense of \cite{CCMT}. The totally disconnected focal hyperbolic groups are also  studied by G. Willis in \cite{WillisSG}, who calls them \textit{scale groups}. They are formed as an ascending HNN extension over a compact open subgroup. In case (ii), the tree need not be locally finite; moreover, the normal subgroup $\cgrp{gKg\inv \mid g \in G}$ is topologically simple by a theorem of R. M\"{o}ller and J. Vonk (\cite[Theorem~2.4]{MollerVonk}).

The following  consequence of Theorem~\ref{intro:cpctend} is immediate, since by definition, a group with a TMS subgroup cannot be compact. 

\begin{corollary}
Let $G$ be a compactly generated \tdlc group without nontrivial compact normal subgroup. If $G$ has a TMS subgroup, then $G$ has infinitely many ends.
\end{corollary}

Proposition~\ref{intro:tree->TMS} may be viewed as a partial converse to Theorem~\ref{intro:cpctend}.  More interestingly, we can express some sufficient conditions to have a TMS subgroup in terms of the local structure of $G$.  In particular, we have the following when $G$ is \defbold{locally of finite quotient type}, meaning that there is a compact open subgroup $U$ of $G$ such that $U$ has only finitely many discrete quotients of each order.

\begin{theorem}\label{intro:cpctend:lqft}
Let $G$ be a nontrivial compactly generated \tdlc group, with $\QZ(G)=\triv$ and with no nontrivial compact normal subgroups, such that $G$ is locally of finite quotient type.  Suppose that there is a nonempty $G$-invariant subset $\mc{Q}$ of $\ldlat(G) \setminus \{0,\infty\}$ such that $G_{\alpha}$ is compact for every $\alpha \in \mc{Q}$ and such that every nonzero element of the subalgebra $\grp{\mc{Q}}$ generated by $\mc{Q}$ lies above some nonzero element of $\mc{Q}$.  Then there is $\alpha \in \mc{Q}$ such that $\rist_G(\alpha) := \bigcap_{\beta \in \grp{\mc{Q}}, \beta \ge \alpha}G_{\beta}$ is a TMS subgroup of $G$, and $G$ falls under case (ii) of Theorem~\ref{intro:cpctend}.
\end{theorem}

The hypothesis that $G_\alpha$ is compact is not a local condition, but seems unavoidable given the known examples.  Specifically, the coloured Neretin groups $G = \mc{N}_{\mathrm{Alt}(2n+5)}$ ($n \ge 0$) constructed by W. Lederle in \cite{Lederle} are in $\ms{S}_{\mc{LD}}$ and are locally of finite quotient type, but they can also be shown to be one-ended by a similar proof to (\cite[Theorem~8.2.6]{Kerg}); in these examples, $G_{\alpha}$ is not compact because it contains isomorphic copies of $G$ itself.  Moreover, $\mc{N}_{\mathrm{Alt}(2n+5)}$ contains as an open subgroup the Burger--Mozes group $U(\mathrm{Alt}(2n+5))^+$, which belongs to $\ms{S}_{\infty}$; so even in the class $\ms{S}$, there is no purely local criterion for the number of ends.

Theorem~\ref{intro:cpctend:lqft}, or rather a special case of its contrapositive (see Corollary~\ref{cor:cpctend:bis}), will be applied in a subsequent article to a large family of groups in $\ms{S}$ (obtained as completions of Kac--Moody groups over finite fields) that were already known to be one-ended, in order to show that they are locally indecomposable.

A known obstacle to weakening the hypotheses of Theorem~\ref{intro:cpctend} is the existence of groups in $\ms{S}$ that are faithful micro-supported, but act geometrically on well-behaved one-ended spaces: see \cite{CapraceRAB}.  However, in the more restrictive context of Theorem~\ref{intro:cpctend:lqft}, where we only consider the decomposition lattice, it would be consistent with known examples to omit the hypothesis that $G$ is locally of finite quotient type and the conclusion that $G$ falls under case (ii) of Theorem~\ref{intro:cpctend}. This leads us to the following. 

\begin{question}
Suppose $G$ and $\mc{Q}$ satisfy the hypotheses of Theorem~\ref{intro:cpctend:lqft}, except that $G$ is not locally of finite quotient type.  Is $\rist_G(\alpha)$ a TMS subgroup of $G$ for some $\alpha \in \mc{Q}$?
\end{question} 

\subsection*{Structure of the article}

In Section~\ref{sec:prelim} we recall some general structure theory of \tdlc groups.  In Section~\ref{sec:arc_geometric} we introduce arc-geometric actions on trees and recall a general construction due to Dicks--Dunwoody, which can be applied in the context of a multi-ended \tdlc group $G$ to obtain an arc-geometric action of $G$ on a tree.  Finally, in Section~\ref{sec:cpctend}, we prove the main results about sources and consequences of TMS subgroups.

\section{Preliminaries}\label{sec:prelim}

\subsection{Local structure theory and dynamics on Stone spaces}

Let $G$ be a \tdlc group.  A subgroup of $G$ is called \defbold{locally normal} if it has open normaliser.  A \defbold{local direct factor} of $G$ is a closed subgroup $K$, such that some open subgroup $O$ of $G$ splits as a topological group as a direct product $O = K \times L$.  Note that every local direct factor is locally normal.  We say $G$ is \defbold{locally indecomposable} if every local direct factor of $G$ is discrete or open.

The \defbold{quasi-centre} $\QZ(G)$ of $G$ consists of all elements of $G$ with open centraliser.  A \tdlc group $G$ is \defbold{[A]-semisimple} if $\QZ(G)=\triv$ and $G$ has no nontrivial abelian locally normal subgroups.

Assume for the moment that $\QZ(G)=\triv$.  Then we define the \defbold{(local) decomposition lattice} $\ldlat(G)$ to be the poset of local direct factors of $G$ ordered by inclusion, modulo the relation $\sim_o$ of \defbold{local equivalence}, where $H \sim_o K$ if $H \cap U = K \cap U$ for some compact open subgroup $U$ of $G$.  The poset $\ldlat(G)$ is then naturally equipped with an action of $G$ by automorphisms, which is induced by the conjugation action of $G$ on its local direct factors.  We define the \defbold{centraliser lattice} $\lcent(G)$ to be the poset of centralisers of locally normal subgroups, modulo local equivalence.

\begin{theorem}[{\cite[Theorem~I]{CRWpart1}}]
Let $G$ be a \tdlc group.
\begin{enumerate}[label=(\roman*)]
\item Suppose $\QZ(G)=\triv$.  Then $\ldlat(G)$ is a Boolean algebra.
\item Suppose $G$ is [A]-semisimple.  Then $\lcent(G)$ is a Boolean algebra containing $\ldlat(G)$ as a subalgebra.
\end{enumerate}
\end{theorem}

Let $\mc{A}$ be a Boolean algebra; write $0$ for the smallest element and $\infty$ for the largest element of $\mc{A}$.  A \defbold{partition} of $\mc{A}$ is a finite set of pairwise disjoint elements of $\mc{A}$ with join $\infty$.  We require all actions on Boolean algebras and topological spaces to be by automorphisms and homeomorphisms respectively.  A Boolean algebra $\mc{A}$ has an associated \defbold{Stone space} $\mf{S}(\mc{A})$, where the points are the ultrafilters $\xi: \mc{A} \rightarrow \{0,1\}$ of $\mc{A}$, and the topology is generated by sets of the form $\{\xi \in \mf{S}(\mc{A}) \mid \xi(\alpha)=1\}$ for $\alpha \in \mc{A}$.  Given an action of a group $G$ on $\mc{A}$, and given $\alpha \in \mc{A}$, we define the \defbold{rigid stabiliser} $\rist_G(\alpha)$ to be the subgroup consisting of all $g \in G$ such that, whenever $\beta \in \mc{A}$ is such that $\beta \ge \alpha$, then $g\beta = \beta$.  Analogously, given an action of a group $G$ on a space $X$, the \defbold{rigid stabiliser} $\rist_G(Y)$ of $Y \subseteq X$ is the pointwise fixator of $X \setminus Y$.  An action of a group $G$ on a Boolean algebra $\mc{A}$, respectively a topological space $X$, is \defbold{micro-supported} if $\rist_G(a)$ acts nontrivially for all $a \in \mc{A} \setminus \{0\}$, respectively for all nonempty open $a \subseteq X$.  In the case of a \tdlc group, we say the action is \defbold{nondiscretely micro-supported} if in addition, the kernel of the action of $\rist_G(a)$ is not open in $\rist_G(a)$, so that every open subgroup of $G$ also has micro-supported action.  An action of a group $G$ on a Boolean algebra $\mc{A}$, respectively a totally disconnected topological space $X$, is \defbold{locally decomposable} if $\prod_{a \in \mc{P}}\rist_G(a)$ is open whenever $\mc{P}$ is a partition of $\mc{A}$, respectively a clopen partition of $X$.

We say a \tdlc group $G$ is \defbold{faithful locally decomposable}, respectively \defbold{faithful micro-supported}, if it has a faithful locally decomposable, respectively faithful nondiscretely micro-supported, action on a Boolean algebra.  We then have a universal $G$-action of this kind, as shown in \cite{CRWpart1}.

\begin{theorem}[See {\cite[Theorem~5.18]{CRWpart1}}]\label{thm:universal_micro}
Let $G$ be a \tdlc group with a compact open subgroup $U$.
\begin{enumerate}[label=(\roman*)]
\item Suppose $G$ has a faithful micro-supported action on some Boolean algebra $\mc{A}$.  Then the following are equivalent:
\begin{enumerate}[label=(\alph*)]
\item $\QZ(G)=\triv$;
\item The action of $U$ is micro-supported, and every nontrivial normal subgroup of $G$ has nontrivial intersection with $U$;
\item $G$ is [A]-semisimple, and $\mc{A}$ is $G$-equivariantly isomorphic to a subalgebra of $\lcent(G)$.  Indeed, the set of rigid stabilisers of the action on $\mc{A}$ forms a subalgebra of the global centraliser lattice of $G$.
\end{enumerate}
\item Up to a $G$-equivariant isomorphism of Boolean algebras, every faithful locally decomposable action of $G$ (if there are any) occurs as the action of $G$ on a subalgebra of $\ldlat(G)$.
\end{enumerate}
\end{theorem}

Given a group $G$ acting on a locally compact topological space $X$, a nonempty subset $Y$ is \defbold{compressible} under the action if for nonempty open $Z$ there exists $g \in G$ such that $gY$ is contained in $Z$.  The action is \defbold{extremely proximal} if every proper compact subspace is compressible.  We make analogous definitions for a group acting on a Boolean algebra via the Stone space.

\subsection{Cayley--Abels graphs, ends and trees}\label{subsec:ends}

A \defbold{graph} $\Gamma$ consists of a vertex set $V\Gamma$, an arc (also known as a directed edge) set $A\Gamma$, origin and terminus functions $o_{\Gamma},t_{\Gamma}:A\Gamma \rightarrow V\Gamma$, and edge reversal $\ol{(-)}_{\Gamma}: A\Gamma \rightarrow A\Gamma$.  (We suppress the subscripts when the graph in question is clear from context.)  In this article all graphs will be simplicial graphs, meaning that $o(a) \neq t(a)$ for all arcs $a$, and the map $a \mapsto (o(a),t(a))$ from $A\Gamma$ to $V\Gamma \times V\Gamma$ is injective; where convenient we will simply identify an arc $a$ with the ordered pair $(o(a),t(a))$, respectively identify the undirected edge $\{a,\ol{a}\}$ with the unordered pair $\{o(a),t(a)\}$.  For the purposes of geometric properties, we also identify a connected graph with its usual geometric realisation; in particular, a connected graph carries a natural metric on its vertices.  A \defbold{bounded} set is one that has finite diameter in the metric.

Given a connected graph $\Gamma$, we can associate compact totally disconnected spaces as follows.  Given $A \subseteq V\Gamma$, let $\delta A$ (or $\delta^\Gamma A$ if the choice of graph $\Gamma$ is not clear from context) be the set of arcs $a \in A\Gamma$ such that $o(a) \not\in A$ and $t(a) \in A$; we say $A$ is \defbold{almost separated} if $\delta A$ is finite.  More generally, in a metric space one has the following coarse geometry concept: $A$ is coarsely almost separated if and only if, for all $d > 0$, the set $\delta_d A$ of points $p$ such that $p \not\in A$ but $d(p,A) \le d$ is bounded.  The two notions of almost separation agree on connected locally finite graphs.

One sees that the set of almost separated subsets of $\Gamma$ is closed under finite intersections, finite unions and complements, so it forms a Boolean algebra $\ms{B} \Gamma$.  (The coarsely almost separated sets also form a Boolean algebra $\ms{A}$, however if $\Gamma$ is not locally finite then $\ms{B} \Gamma$ can be properly contained in $\ms{A}$.)  If $\Gamma$ itself is unbounded, the bounded almost separated subsets form an ideal of $\ms{B}\Gamma$; write $\ol{\ms{B}\Gamma}$ for the quotient of $\ms{B}\Gamma$ by this ideal.  We then have associated Stone spaces $\ol{\Gamma} = \mf{S}(\ms{B} \Gamma)$ and $\ol{\Ends} \Gamma = \mf{S}(\ol{\ms{B} \Gamma})$; equivalently $\ol{\Ends} \Gamma$ can be regarded as the closed subspace of $\mf{S}(\ms{B} \Gamma)$ consisting of the ultrafilters that are zero on bounded subsets.  (If $\Gamma$ itself is bounded, we can take $\ol{\Ends} \Gamma := \emptyset$.)  Given $\xi \in \ol{\Ends} \Gamma$ and $A \in \ms{B}\Gamma$, we will say ``$\xi$ is in $A$'' or ``$A$ is in $\xi$'' interchangeably to mean $\xi(A)=1$.

Using coarse almost separation, one sees that if $\Gamma$ is locally finite then $\ol{\Ends} \Gamma$ is a coarse geometric invariant of the graph, in other words it is preserved by passing to a cobounded subset with a large-scale equivalent metric.  In addition, the space $\ol{\Ends} \Gamma$ is naturally homeomorphic to the usual notion of the space of ends $\Ends \Gamma$ of a geodesic metric space, as defined by equivalence classes of rays (see for instance \cite[Chapter~I.8, Proposition~8.29]{BHCAT0}).

A \defbold{tree} is a simply connected graph; a \defbold{subtree} is a nonempty subgraph that is also a tree.  Given a tree $T$ and an arc $e$ of $T$, we define the associated \defbold{half-tree} $T_e$ to be the induced graph on the vertices $v$ of the tree such that $d(t(e),v) < d(o(e),v)$; in particular, note that $T_e$ is an almost separated set, with $\delta T_e = \{e\}$.  It is sometimes useful to distinguish within $\ol{\Ends} T$ the subspace of \defbold{(geometric) ends} $\Ends T$ of $T$, consisting of those $\xi \in \ol{\Ends} T$ containing an infinite descending sequence of half-trees; these correspond in a natural way to equivalence classes of geodesic rays in the tree.

Now let $G$ be a compactly generated \tdlc group.  An action of $G$ on a metric space $X$ is called \defbold{geometric} if the following conditions are satisfied:
\begin{enumerate}[label=(\alph*)]
\item The action is \textbf{isometric}: $G$ acts by isometries;
\item The action is \textbf{proper}: For all $n < \infty$ and $x \in X$, the set of $g \in G$ such that $d(x,gx) < n$ is a neighbourhood of the identity with compact closure;
\item The action is \textbf{cobounded} (the term \textbf{cocompact} is used if $X$ is locally compact): There is a distance $n$ such that for all $x,y \in X$, there is $g \in G$ such that $d(x,gy) < n$.
\end{enumerate}
A metric space equipped with a geometric action of $G$ is called a \defbold{$G$-metric space}, and a \defbold{$G$-metric} on a $G$-set is a metric with respect to which the action of $G$ is geometric.  All proper geodesic $G$-metric spaces have the same quasi-isometry type, which is also the quasi-isometry type of $G$ as defined intrinsically; indeed, if $G$ acts geometrically on the proper geodesic metric space $M$ then $g \mapsto gm$ is a quasi-isometry for any $m \in M$.  (See \cite[Sections 4.B and 4.C]{CdH16}.)  In particular we can define the \defbold{space of ends $\Ends G$ of $G$} to be $\ol{\Ends} M$, where $M$ is any proper geodesic $G$-metric space.  There is then a natural induced action of $G$ on $\Ends G$; from the construction it is clear that $G$ acts by homeomorphisms.

With \tdlc groups, it is useful to think about geometric actions more combinatorially.  A \defbold{$G$-metric graph} is a connected locally finite graph equipped with a geometric action of $G$ with respect to the graph metric, and a \defbold{geometric $G$-set} is a set $X$ equipped with a permutation action of $G$, such that point stabilisers are compact open and $G$ has finitely many orbits.  A \defbold{Cayley--Abels graph} for $G$ is a connected locally finite simple graph $\Gamma$ equipped with an action of $G$ by isometries that is vertex-transitive, with compact open vertex stabilisers.  In particular, every Cayley--Abels graph is a $G$-metric graph.

\begin{lemma}[{See e.g. \cite[Section~4.1]{CRWpart2}}]\label{lem:CayleyAbels}
Let $G$ be a compactly generated \tdlc group, let $U$ be a compact open subgroup of $G$ and let $A$ be a compact symmetric subset of $G$ such that $G = \grp{U,A}$.  Then there is a finite symmetric subset $B$ of $G$ such that
$$ BU = UB = UBU = UAU.$$
Moreover, for any such subset $B$, we have $G = \grp{B}U$ and the coset space $G/U$ is the set of vertices of a Cayley--Abels graph $\Gamma$ for $G$, where $gU$ is adjacent to $hU$ in $\Gamma$ if and only if $gU \neq hU$ and $Uh\inv gU \subseteq UBU$.
\end{lemma}

Any transitive geometric $G$-set is of the form $G/U$ for some compact open subgroup $U$ of $G$; Lemma~\ref{lem:CayleyAbels} then supplies a Cayley--Abels graph $\Gamma$ with $V\Gamma = G/U$.  It is then straightforward to see that given any geometric $G$-set $X$, there is a geometric $G$-graph $\Gamma$ with $V\Gamma = X$ (restricting to a Cayley--Abels graph on each $G$-orbit of vertices), yielding a $G$-metric on $X$.  In particular, the space of ends $\Ends G$ of any compactly generated \tdlc group $G$ is realised as a $G$-space by $\Ends \Gamma$ for a Cayley--Abels graph $\Gamma$ of $G$.

\subsection{Local finiteness properties}

\begin{definition}
Say a profinite group $U$ is of \defbold{finite quotient type} if for each natural number $n$, there are only finitely many open normal subgroups of $U$ of index $n$.  Say a \tdlc group $G$ is \defbold{locally of finite quotient type} if every compact open subgroup of $G$ is of finite quotient type.
\end{definition}

To determine if $G$ is locally of finite quotient type, it suffices to consider a single compact open subgroup.

\begin{lemma}\label{lem:fqt_local}
Let $G$ be a \tdlc group.  Suppose there is some compact open subgroup $U$ of $G$ of finite quotient type.  Then for every compact open subgroup $V$ of $G$ and $n \in \bN$, there are only finitely many open subgroups of $V$ of index $n$.  In particular, $G$ is locally of finite quotient type.
\end{lemma}

\begin{proof}
Every open subgroup of $U$ of index $n$ contains an open normal subgroup of index dividing $n!$.  Consequently there are only finitely many open subgroups of $U$ of index $n$.

Now consider an arbitrary compact open subgroup $V$ of $G$.  Given an open subgroup $W$ of $V$ of index $n$, we see that $W$ contains an open subgroup of $U$ of index at most $n|U:U \cap V|$.  Thus for a fixed $n$, there are only finitely many such subgroups $W$ of $V$.
\end{proof}

We note a property of groups that are locally of finite quotient type that will be useful later.

\begin{lemma}\label{lem:strong_fin_decomp}
Let $G$ be a nontrivial compactly generated \tdlc group that is locally of finite quotient type with $\QZ(G)=\triv$ and such that $\bigcap_{g \in G}gOg\inv = \triv$ for some open subgroup $O \le G$.  Then each compact open subgroup $U$ of $G$ has only finitely many direct factors.
\end{lemma}

\begin{proof}
Let $U$ be a compact open subgroup of $G$.  Since $\bigcap_{g \in G}gOg\inv = \triv$, it follows by \cite[Proposition~4.6]{CRWpart2} that the composition factors of $U$ are of bounded order; let $k$ be the largest order that occurs.  Then $U$ has only finitely many quotients of order at most $k$, and hence only finitely many simple quotients.  In any profinite group $P$, a standard argument shows that any proper closed normal subgroup is contained in an open normal subgroup $Q$ such that $P/Q$ is simple; by considering the simple quotients of $U$, we see that $U$ cannot be written as an infinite direct product of nontrivial profinite groups.  By \cite[Proposition~4.11]{CRWpart1} it follows that $U$ has only finitely many direct factors.
\end{proof}

\begin{remark}
A profinite group $U$ is \defbold{topologically finitely generated} if it has a dense finitely generated subgroup; it is then easy to deduce that $U$ has only finitely many open subgroups of each index.  Given a pro-$p$ group, having finite quotient type (indeed, having only finitely many open normal subgroups of index $p$) is equivalent to being topologically finitely generated, since the quotient by the Frattini subgroup is elementary abelian.  For general profinite groups on the other hand, finite quotient type is strictly more general than topological finite generation.  For example, every just infinite profinite group is of finite quotient type (see for instance \cite[Corollary~2.5]{ReiJI}), whereas there are examples (see for instance \cite{WilHJI}) of (hereditarily) just infinite profinite groups that are not topologically finitely generated.
\end{remark}

\section{Arc-geometric actions of \tdlc groups on trees}\label{sec:arc_geometric}

In this section we analyse the space of ends of a compactly generated \tdlc group via actions on trees.  This approach is well known, and finds its origin in the analogue of Stallings' theorem in this context, which is due to H.~Abels \cite{Abels-74}. Moreover,  an analysis along similar lines can also be found in the article \cite[\S 3]{KronMoller} of Kr\"{o}n--M\"{o}ller, which was in turn inspired by the Dicks--Dunwoody approach (see \cite{DDbook}) to the space of ends of an abstract group.  However, the authors of the present article found it useful to give a new presentation of the ideas, in particular in order to phrase the results in terms of arc-geometric actions and to highlight the dynamics of the action on $\Ends G$.

A group acting on a tree, preserving no proper subtree, has one of a few possible structures.  The next proposition is already known in some form (indeed similar results are known in greater generality, see for example \cite[\S 3.A--3.B]{CCMT}), but we could not find a clear reference for the result as stated here, so we include a proof for clarity.

\begin{prop}\label{prop:tree_dense}
Let $T$ be a tree with more than two vertices and let $G \le \Aut(T)$.  Suppose that $G$ preserves no proper subtree of $T$.  Then exactly one of the following holds.
\begin{enumerate}[label=(\roman*)]
\item $T$ is a line and $G$ is either cyclic or infinite dihedral, acting geometrically on $T$.
\item $G$ fixes exactly one end $\xi \in \Ends T$ and has the form $P \rtimes \grp{l}$, where $l$ is a translation towards $\xi$ and the $P$-orbits on $VT$ are exactly the horospheres around $\xi$.  Consequently, the action of $G$ on $\ol{\Ends} T \setminus \{\xi\}$ is extremely proximal.
\item $G$ has no fixed points in $\Ends T$, and for all pairs $(a,b)$ of arcs of $T$ there is $g \in G$ such that $T_{ga} \subseteq T_b$.  Consequently, the action of $G$ on $\Ends T$ is faithful and the action on $\ol{\Ends} T$ is extremely proximal.
\end{enumerate}
Moreover, in cases (ii) and (iii) then $\ol{\Ends} T$ is infinite and perfect (that is, it has no isolated points) and $G$ acts faithfully on $\Ends T$; moreover, every nontrivial normal subgroup of $G$ has unbounded orbits on $T$.
\end{prop}

\begin{proof}
If $T$ is a line, then the fact that $G$ preserves no proper subtree means that $G$ contains a translation.  It is then clear that $G$ is a cyclic or infinite dihedral group acting geometrically, as in (i).  We may now assume $T$ is not a line; the hypotheses ensure $T$ is leafless, so $T$ has more than two ends.  Since $T$ is leafless, it is also easy to see that every half-tree belongs to a geometric end of $T$, so $\Ends T$ is dense in $\ol{\Ends} T$.

Suppose $G$ fixes an end $\xi$; if $G$ fixed another end $\xi'$ then $G$ would fix the line between them, which is a proper subtree of $T$.  Thus $\xi$ is unique.  Let $\beta: VT \rightarrow \bZ$ be a function such that if $a$ is an arc pointing towards $\xi$ then $\beta(o(a)) - \beta(t(a)) =1$.  By considering vertices along two geodesic rays from some initial vertex $v$, where one ray represents $\xi$ and the other represents an end other than $\xi$, we see that $\beta$ is surjective.  Let $r: [0,\infty) \rightarrow T$ be a geodesic ray representing $\xi$.  Since $G$ fixes $\xi$, given $g \in G$ there exist $n_0$ and $t$ such that $gr(n) = r(n+t)$ for all $n \ge n_0$.  From the graph structure we have $t \in \bZ$, and then one sees that in fact $\beta(v)-\beta(gv) = t$ for all $v \in VT$.  Thus $G = P \rtimes \grp{l}$ where $l$ is either trivial or a translation towards $\xi$ and $P$ is the setwise stabiliser of each of the horospheres $\beta\inv(n)$ around $\xi$.  We can moreover rule out the case that $l$ is trivial by noting that $P$ preserves a proper subtree, namely the subtree spanned by $\beta\inv(-\bN)$.  Now let $n \in \bZ$ and suppose that $P$ acts intransitively on $\beta\inv(n)$, say $X_n$ is a proper nonempty $P$-invariant subset.  We see that $l$ translates along an axis $L$, and then the union of the $P$-orbit of $L$ is a $G$-invariant subtree $T'$ such that $VT' \cap \beta\inv(n)$ is contained in either $X_n$ or its complement in $\beta\inv(n)$.  In either case we have a proper $G$-invariant subtree, contradicting our hypothesis.  So in fact $P$ acts transitively on $\beta\inv(n)$ for each $n \in \bZ$.

Let $Y$ and $Z$ be nonempty compact open subspaces of $\ol{\Ends} T \setminus \{\xi\}$.  We see that $Y \subseteq \ol{T_a}$ and $\ol{T_b} \subseteq Z$ for some half-trees $T_a$ and $T_b$ of $T$, such that $\xi \not\in \Ends T_a$.  From the fact that $P$ acts transitively on every horosphere, we see that it acts minimally on $\ol{\Ends} T \setminus \{\xi\}$.  In particular, writing $\xi'$ for the end of $L$ other than $\xi$, we see that there is $g \in P$ such that $g\xi' \in \Ends T_b$.  We then have $gl^{-n}g\inv T_a \subseteq T_b$, and hence $gl^{-n}g\inv Y \subseteq Z$, for some $n \ge 0$.  Thus the action on $\ol{\Ends} T \setminus \{\xi\}$ is extremely proximal.

The remaining case is that $G$ fixes no ends of $T$.  Let $a$ and $b$ be arcs of $T$.  Since $G$ does not fix an end, we see that there is $k \in G$ such that $a$ and $ka$ point away from each other, that is, $T_a$ and $T_{ka}$ are disjoint.  Now let $H$ be the set of $h \in G$ such that $o(ha) \in T_b$ and $ha \not\in \{b,\ol{b}\}$.  Then $H \neq \emptyset$ by \cite[Lemme~4.1]{Tits70}; taking $h \in H$, if $T_{ha}$ is not contained in $T_b$, then $T_{\ol{ha}} \subseteq T_b$ and hence $T_{hka} \subseteq T_b$.  Thus there is $g \in G$ such that $T_{ga} \subseteq T_b$.

To show the action is extremely proximal it suffices to show, for any two proper nonempty clopen subsets $Y$ and $Z$ of $\ol{\Ends} T$, that there is $g \in G$ such that $gY \subseteq Z$.  As before we have $Y \subseteq \ol{T_a}$ and $\ol{T_b} \subseteq Z$ for some half-trees $T_a$ and $T_b$ of $T$.  By the previous paragraph, there is then $g \in G$ such that $T_{ga} \subseteq T_b$, from which it follows that $gY \subseteq Z$.  We have now shown that the action of $G$ on $\ol{\Ends} T$ is extremely proximal, proving (iii).

Now suppose we are in case (ii) or (iii); let $B$ be the set of vertices of $T$ of degree at least $3$.  In either case it is clear that $\ol{\Ends} T$ is infinite.  In case (ii), we see that some horosphere, say $\beta\inv(0)$, intersects $B$; it then follows that $\beta\inv(mn) \subseteq B$ for all $n \in \bZ$, where $m$ is the translation length of $l$, so $\ol{\Ends} T$ is perfect.  In case (iii), $\ol{\Ends} T$ is perfect because it is a compact minimal $G$-space, so again the convex hull of $B$ is the whole tree.  Now let $g \in G$ act trivially on $\Ends T$.  Then given $v \in B$, there are ends $\xi_1,\xi_2,\xi_3$ forming the corners of an ideal triangle of $T$, such that $v$ is the unique point lying on all three sides of the triangle; in particular, $g$ fixes $B$ pointwise.  Since the convex hull of $B$ is $T$, we deduce that $g=1$, as required.  Given a normal subgroup $N$ with bounded orbits on $T$, we see that $N$ fixes a vertex or inverts an edge; the fact that $G$ preserves no subtree then implies that $N$ is trivial.
\end{proof}

Given a compactly generated \tdlc group $G$ acting on a leafless tree $T$, we say the action is \defbold{arc-geometric} if the arcs form a geometric $G$-set.  This situation can be characterised as follows:

\begin{lemma}[{See also \cite[Theorem~11]{KronMoller}}]\label{lem:tree_geometric}
Let $G$ be a compactly generated \tdlc group acting on a leafless tree $T$.
\begin{enumerate}[label=(\roman*)]
\item The action is arc-geometric if and only if $G$ has compact open arc stabilisers and preserves no proper subtree.
\item Suppose the action of $G$ on $T$ is arc-geometric.  Then there is a $G$-metric graph $\Gamma$ with vertex set $V\Gamma = AT$ where for all arcs $(a,b)$ of $\Gamma$, either $b = \ol{a}$ or $o_T(a) = o_T(b)$.  Moreover, for almost separated set $X$ in $VT$, the set $A_TX$ of all arcs of $T$ between vertices of $X$ is almost separated in $\Gamma$, so we have a $G$-equivariant quotient map $\Ends G \rightarrow \ol{\Ends} T$.
\end{enumerate}
\end{lemma}

\begin{proof}
Suppose the action is arc-geometric; by definition, $G$ has compact open arc stabilisers.  Now let $T'$ be a proper $G$-invariant subtree.  Then since $T$ is leafless, there is some geodesic ray $r$ of $T$ such that $r(0)$ is a vertex of $T'$, but thereafter $r$ is outside $T'$.  In particular, we see that the distance from $r(n)$ to $T'$ tends to infinity, meaning that the arcs in the image of $r$ lie in infinitely many $G$-orbits, a contradiction.  Thus $G$ preserves no proper subtree.

Conversely, suppose that $G$ has compact open arc stabilisers and preserves no proper subtree of $T$.  After taking a quotient with compact kernel we may assume $G$ acts faithfully on $T$.  By \cite[Proposition~6.6]{ReidSmith} one sees that $G$ has finitely many orbits on $AT$, and $G_v$ is compactly generated for each vertex $v$.  Thus $G$ has arc-geometric action on $T$.  Consider now the action of $G_v$ on $o\inv(v)$: given $a,b \in o\inv(v)$, then $a$ and $b$ are in the same $G$-orbit if and only if they are in the same $G_v$-orbit, and $G_a = (G_v)_a$.  Thus $o\inv(v)$ is a geometric $G_v$-set.  Now take $v_1,\dots,v_k$ to be representatives of the $G$-orbits on $VT$, and for $1 \le i \le k$ make a connected locally finite graph $\Gamma_{v_i}$ with vertex set $o\inv(v_i)$, such that $G_{v_i}$ acts geometrically on $\Gamma_{v_i}$.  We then define a graph $\Gamma$ with vertex set $AT$ and arcs
\[
A\Gamma = \bigcup_{a \in AT}(a,\ol{a}) \cup \bigcup^k_{i=1}\bigcup_{g \in G}gA\Gamma_{v_i}.
\]
It is then easy to see that $\Gamma$ is connected, locally finite and $G$-invariant, so that $G$ acts on $\Gamma$ geometrically.  By construction, if $(a,b)$ is an arc of $\Gamma$ then either $b = \ol{a}$ or $o_T(a) = o_T(b)$.

Let $X$ be an almost separated set of vertices in $T$ and let $A_TX$ be the set of arcs $a \in AT$ such that $\{o_T(a),t_T(a)\} \subseteq X$.  Consider an arc $(a,b)$ of $\Gamma$ such that $a \not\in A_TX$ and $b \in A_TX$.  Then $b \neq \ol{a}$, so $o_T(a) = o_T(b)$.  However, the only way $a$ can be outside $A_TX$ while having the same origin as an arc in $A_TX$ is if $o_T(a) \in X$ and $t_T(a) \not\in X$, in other words, $a \in \delta^T X$.  Then by the fact that $\Gamma$ is locally finite, there are only finitely many possibilities for the pair $(a,b)$, so $\delta^\Gamma X^*$ is finite.  Thus $A_TX$ is almost separated as a subset of $V\Gamma$.  We deduce that 
\[
\{A_TX \mid X \subseteq VT \text{ almost separated}\}
\]
is a subset $\ms{A}$ of $\ms{B}\Gamma$.  Note also that if $X$ is bounded, then actually it is finite and consists of vertices of finite degree, so $A_TX$ is also finite.  Conversely, if $AX$ is bounded then it is finite; since $AX$ includes all but finitely many of the $T$-arcs incident with $X$, we deduce that $X$ is finite.  Moreover, we find that for any almost separated sets $X_1,X_2 \subseteq VT$, then $A_T(X_1 \cup X_2)$ has finite difference with $A_TX_1 \cup A_TX_2$, and similarly for intersections and complements.  Thus the map $X \mapsto A_TX$ induces an injective homomorphism $\iota: \ol{\ms{B}T} \rightarrow \ol{\ms{B}\Gamma}$; the dual map $\hat{\iota}$ is then a quotient map from $\ol{\Ends} \Gamma$ (which can be identified with $\Ends G$) to $\ol{\Ends} T$.
\end{proof}

In particular, we see that if the compactly generated \tdlc group $G$ has an arc-geometric action on a leafless tree $T$, then it has at least $|\ol{\Ends} T|$ ends.  Conversely, if $G$ has more than one end, then the space of ends can be approximated using arc-geometric actions.  It is useful to divide into three cases:

\begin{enumerate}[label=(\Alph*)]
\item $G$ acts geometrically on a line.
\item $G$ has a vertex-transitive action on a locally finite tree $T$ with compact open stabilisers, fixing exactly one end.
\item $G$ has more than end but does not satisfy (A) or (B).
\end{enumerate}

In case (A) it is clear that $G$ is compact-by-$D$ where $D$ is either infinite cyclic or infinite dihedral, and that $|\Ends G|=2$.  In particular, Lemma~\ref{lem:tree_geometric} ensures that (A) and (B) are mutually exclusive.

For case (B), say a compactly generated \tdlc group $G$ is a \defbold{scale group} if $G$ admits a faithful vertex-transitive action on a locally finite tree $T$ with compact open stabilisers, fixing exactly one end.  We say $G$ is an \defbold{almost scale group} if $G/N$ is a scale group for some compact normal subgroup $N$.  Every almost scale group is a focal hyperbolic group in the sense of \cite{CCMT}. The results from loc. cit. imply the following characterisation.

\begin{prop}\label{prop:scale_group}
Let $G$ be a compactly generated \tdlc group.  Then $G$ is an almost scale group if and only if it has an arc-geometric action on a leafless tree $T$, fixing an end, such that $T$ is not a line.  If $G$ has such an action, then $T$ is locally finite and $\Ends G \cong \Ends T$ as $G$-spaces.

Moreover, if $G/N$ is a scale group for some compact normal subgroup $N$, then $G$ has exactly one fixed point in $\Ends G$ and $N$ is the largest compact normal subgroup of $G$.
\end{prop}

\begin{proof}
If $G/N$ is a scale group for some compact normal subgroup $N$, then $G$ clearly has the specified action on a locally finite tree and also on its space of ends.  Since $G/N$ acts vertex-transitively on an infinitely ended tree, its unique compact normal subgroup is trivial, so $N$ is the largest compact normal subgroup of $G$.

Conversely, suppose $G$ admits an arc-geometric action on a leafless tree $T$, fixing an end $\xi$, such that $T$ is not a line.  Then by Lemma~\ref{lem:tree_geometric} the action preserves no proper subtree, so we are in case (ii) of Proposition~\ref{prop:tree_dense}, with $G = P \rtimes \grp{l}$.  If $a$ is an arc pointing towards the fixed end $\xi$, then $G_a = G_{o(a)}$.  Moreover, the fact that $P$ acts transitively on each horosphere implies that $G_{o(a)}$ must act transitively on $o\inv(o(a)) \setminus \{a\}$.  Since $G_{o(a)}$ is compact and arc stabilisers are open, we deduce that $o\inv(o(a))$ is finite, and thus the tree is locally finite and hence a $G$-metric graph, ensuring that $\Ends G \cong \Ends T$ as $G$-spaces.  In particular, $G$ is hyperbolic, hence focal hyperbolic in the sense of \cite{CCMT}. It follows from \cite[Lemma~5.1]{CCMT} that $G$ has a largest compact normal subgroup $N$, and then it follows from \cite[Theorem~7.1(b) and Theorem~7.3]{CCMT} that $G/N$ is a scale group.
\end{proof}

It remains to give the approximation of $\Ends G$ using arc-geometric actions on trees in case (C).  Say that a compactly generated \tdlc group $G$ is of \defbold{general infinitely-ended type} if $G$ has more than one end, does not act geometrically on a line, and is not an almost scale group.  In this case, we can appeal to a construction by Dicks--Dunwoody.  Given a graph $\Gamma$ and $A \subseteq V\Gamma$, let $\delta A$ be the set of arcs $a$ such that $o(a) \not\in A$ and $t(a) \in A$; notice that $\delta A$ is finite if and only if $A$ is almost separated.  Write $\ms{B}_n\Gamma$ for the subalgebra of $\ms{B}\Gamma$ generated by the almost separated sets $A$ such that $|\delta A| \le n$; note that $\ms{B}_n\Gamma$ is $G$-invariant and that $(\ms{B}_n\Gamma)_{n \ge 0}$ is an ascending sequence of subalgebras with union $\ms{B}\Gamma$.  By \cite[Theorem~II.2.20]{DDbook}, there is an ascending sequence $(\mc{R}_n)_{n \ge 1}$ of nonempty $G$-invariant subsets of $\ms{B}\Gamma$, such that $\mc{R}_n$ generates $\ms{B}_n\Gamma$ as a subalgebra and has the poset structure of the set of half-trees of a tree: specifically, it carries an order-reversing involution $c$ (which is just complementation in $V\Gamma$), such that for all $A,B \in \mc{R}_n$ exactly one of the following holds:
\[
A < B,\  A = B, \  A > B, \ A < B^c, \ A = B^c, \ A > B^c,
\]
and such that for all $A,B \in \mc{R}_n$ there are only finitely many $C \in \mc{R}_n$ such that $A < C < B$.  The elements $A \in \mc{R}_n$ can also be chosen so that both $A$ and $A^c$ span connected subgraphs of $V\Gamma$.  There is then an associated tree $T^{(n)} = T(\mc{R}_n)$ such that the half-trees correspond to the poset $\mc{R}_n$ with the same partial order, yielding a $G$-equivariant embedding $\iota_n: \ms{B}T^{(n)} \rightarrow \ms{B}\Gamma$ with image $\ms{B}_n\Gamma$.  The tree $T^{(n)}$ is also equipped with a canonical $G$-equivariant map $\varphi_n: V\Gamma \rightarrow VT^{(n)}$, with the property that $\varphi_n(v)$ belongs to the half-tree corresponding to $A \in \mc{R}_n$ if and only if $v \in A$.  In particular, note that $\varphi_n(v) \neq \varphi_n(w)$ if and only if $v$ and $w$ are separated by some element of $\ms{B}_n\Gamma$.  (See for instance \cite[\S 3.1.3]{KronMoller}.)

\begin{lemma}\label{lem:DD_subtree}
Let $G$ be a compactly generated \tdlc group, let $\Gamma$ be a Cayley--Abels graph for $G$ and let the tree $T^{(n)}$ be as above for some $n \ge 1$.  Then the action of $G$ on $T^{(n)}$ preserves no proper subtree.  If $T^{(n)}$ has more than one vertex, then the action of $G$ on $T^{(n)}$ is arc-geometric.
\end{lemma}

\begin{proof}
We may assume $T = T^{(n)}$ has more than one vertex.  We observe that $\varphi_n(V\Gamma)$ intersects every half-tree of $T$, since every element of $\mc{R}_n$ contains a vertex.  Since $G$ acts transitively on $V\Gamma$, it also acts transitively on $\varphi_n(V\Gamma)$; it follows by \cite[Lemme~4.1]{Tits70} that $G$ preserves no proper subtree.

Given $a \in AT$, then $G_a$ is the setwise stabiliser of the element $R_a \in \mc{R}_n$ corresponding to the half-tree $T_a$.  Since $\delta R_a$ is bounded, $G_a$ has bounded orbits on $\Gamma$; since $G$ acts geometrically on $\Gamma$ it follows that $G_a$ is compact.  On the other hand, since $\delta R_a$ consists of finitely many arcs, $G_a$ is open.  Thus by Lemma~\ref{lem:tree_geometric}, the action is arc-geometric.
\end{proof}

We thus obtain an approximation of $\Ends G$ by boundaries of trees on which $G$ acts arc-geometrically.

\begin{prop}\label{prop:DD_approximation}
Let $G$ be a compactly generated \tdlc group of general infinitely-ended type and let $\Gamma$ be a Cayley--Abels graph for $G$.  Then $G$ has a unique largest compact normal subgroup $N$, and the action of $G$ on $\Ends \Gamma$ is extremely proximal with kernel $N$.  Moreover, constructing the trees $T^{(n)}$ as above, then there is some $t$ such that for all $n \ge t$, the dual of the inclusion map $\iota_n: \ms{B}T^{(n)} \rightarrow \ms{B}\Gamma$ restricts to a quotient map $\widehat{\iota}^*_n: \Ends \Gamma \rightarrow \ol{\Ends} T^{(n)}$.  In particular, $\Ends G \cong \varprojlim_{n \ge t} \ol{\Ends} T^{(n)}$ as $G$-spaces.
 \end{prop}

\begin{proof}
Consider two proper nonempty clopen subsets $Y$ and $Z$ of $\Ends \Gamma$.  Via Stone duality, $Y$ and $Z$ correspond to elements $\ol{A}_Y$ and $\ol{A}_Z$ of $\ol{\ms{B}\Gamma} \setminus \{0,\infty\}$.  Choose $n$ large enough that $T = T^{(n)}$ has more than one vertex and $\ol{A}_Y$ and $\ol{A}_Z$ belong to the image of $\ms{B}_n\Gamma$ in $\ol{\ms{B}\Gamma}$, and take representatives $A_Y$ and $A_Z$ of $\ol{A}_Y$ and $\ol{A}_Z$ respectively in $\ms{B}_n\Gamma$.  Then $\varphi_n(A_Y)$ and $\varphi_n(A_Z)$ are each described as subsets of $\varphi_n(V\Gamma)$ by taking unions of intersections with finitely many half-trees of $T$, and hence there are half-trees $T_a$ and $T_b$ of $T$ such that
\[
\varphi_n(A_Y) \subseteq T_a \; \text{ and } \; T_b \cap \varphi(V\Gamma) \subseteq T_b.
\]
By Proposition~\ref{prop:tree_dense} we see that $gT_a \subseteq T_b$ for some $g \in G$, which then implies $gA_Y \subseteq A_Z$ (since $A_Y$ and $A_Z$ are saturated with respect to $\varphi_n$).  Thus $G$ has extremely proximal action on $\Ends \Gamma$.  We also see by Proposition~\ref{prop:tree_dense} that if $T^{(n)}$ has more than one vertex, then the kernel $N$ of the action of $G$ on $T^{(n)}$ is the largest compact normal subgroup of $G$.

Now let $t$ be large enough that $T^{(n)}$ has more than one vertex for all $n \ge t$, and consider the image of $\Ends \Gamma$ under the dual map $\widehat{\iota}_n$ from $\ol{\Gamma}$ to $\ol{T^{(n)}}$ for $n \ge t$.  We see that $G$ has extremely proximal, in particular, minimal, action on $\widehat{\iota}_n(\Ends \Gamma)$ and on $\ol{\Ends} T^{(n)}$; thus the two subspaces of $\ol{T^{(n)}}$ are either equal or disjoint.  Moreover, taking $g \in G$ with hyperbolic action on $T^{(n)}$, then by considering the action of $g$ on $\varphi_n(V\Gamma)$, we see that $g$ has an attracting end $\xi$ on $\Gamma$, which is mapped to the attracting end of $g$ on $T^{(n)}$.  Thus $\widehat{\iota}_n(\Ends \Gamma)$ and $\ol{\Ends} T^{(n)}$ are not disjoint, so they are equal.  Since the spaces are compact Hausdorff and $\widehat{\iota}_n$ is continuous, it restricts to a quotient map $\widehat{\iota}^*_n: \Ends \Gamma \rightarrow \ol{\Ends} T^{(n)}$.  In particular, identifying $\Ends \Gamma$ with $\Ends G$, we recover $\Ends G$ as the desired inverse limit, and the kernel of the action is again $N$.
\end{proof}

To summarise this section, here are the possible large-scale structures of a compactly generated \tdlc group $G$, regarded from the perspective of the space of ends.

\begin{prop}\label{prop:end_types}
Let $G$ be a compactly generated \tdlc group and let $N$ be the kernel of the action of $G$ on $\Ends G$.  Then exactly one of the following holds:
\begin{description}
\item[(Compact)] $G$ is compact and $G=N$;
\item[(One-ended)] $|\Ends G| = 1$ and $G=N$;
\item[(Two-ended)] $|\Ends G| = 2$ and $G/N \in \{\bZ, \bZ/2\bZ \ast \bZ/2\bZ\}$, which implies in particular that $G$ acts geometrically on a line;
\item[(Focal hyperbolic)] $\Ends G$ is homeomorphic to the Cantor set, $N$ is compact and $G/N$ is a scale group, fixing exactly one point $\xi$ in $\Ends G$ and acting faithfully and extremely proximally on $\Ends G \setminus \{\xi\}$;
\item[(General infinitely-ended)]  We have a sequence of trees $T_i$, such that on each tree, $G$ acts arc-geometrically on $T_i$ and extremely proximally with compact kernel $N$ on $\ol{\Ends} T_i$, and $\Ends G \cong \varprojlim \ol{\Ends} T_i$ as $G$-spaces.
\end{description}
\end{prop}

We remark that if $G$ is noncompact with no infinite cyclic or infinite dihedral quotient (for example, $G \in \ms{S}$), then only the one-ended and general infinitely-ended cases can occur.

\section{TMS subgroups}\label{sec:cpctend}

\subsection{A sufficient condition for a micro-supported action on a tree}\label{subsec:cpctend}

In this subsection we establish a sufficient condition for a compactly generated \tdlc group $G$ to have an arc-geometric action on a tree with micro-supported action on the boundary.  In order to state the theorem, we define a certain kind of subgroup of a \tdlc group.

\begin{definition}
Let $G$ be a \tdlc group, let $K$ be a subgroup of $G$ and let $U$ be a compact open subgroup of $G$.  We say $(K,U)$ is a \defbold{TMS pair} of $G$ if it has the following properties:
\begin{enumerate}[label=(\alph*)]
\item For each compact normal subgroup $N$ of $G$, and each $G$-conjugate $V$ of $U$, the set of $g \in G$ such that $gKg\inv \le UN$ but $gKg\inv \nleq V$ has compact closure.
\item The set of $g \in G$ such that $gKg\inv \le U$ does not have compact closure.
\item For each compact normal subgroup $N$ of $G$, the index $|K:K \cap N|$ is infinite.
\end{enumerate}
We say $K \le G$ is a \defbold{TMS subgroup} of $G$ if it forms a TMS pair $(K,U)$ with some compact open subgroup $U$ of $G$.
\end{definition}

Here is the main theorem of this section.

\begin{theorem}\label{cpctend}
Let $G$ be a compactly generated \tdlc group with a compact normal subgroup $M$, such that there is a TMS subgroup $K/M$ of $G/M$.  Then $G$ has a unique largest compact normal subgroup $N$ and $G/N$ is [A]-semisimple, with faithful action on $\Ends G \cong \Ends G/N$ such that every nonempty open set contains the support of some $G/N$-conjugate of $K/N$.  Moreover, one of the following holds.
\begin{enumerate}[label=(\roman*)]
\item $G/N$ is a scale group.
\item We have a sequence of trees $(T_i)_{i \ge 0}$ on which $G$ acts, with the following properties:
\begin{enumerate}[label=(\alph*)]
\item For each $i$, the action of $G$ on $T_i$ is arc-geometric, with kernel $N$ and with the support of $K$ confined to some half-tree;
\item As $G$-spaces, $\Ends G \cong \varprojlim \ol{\Ends} T_i$ and the action on $\Ends G$ is extremely proximal.
\end{enumerate}
\end{enumerate}
\end{theorem}

We begin the proof with a series of lemmas, starting with some general observations on the TMS property.

\begin{lemma}\label{lem:tms_technicalities}
Let $G$ be a \tdlc group and let $(K,U)$ be a TMS pair for $G$.
\begin{enumerate}[label=(\roman*)]
\item $K$ is infinite with compact closure, so in particular it is nondiscrete, and $(g\ol{K}g\inv,U)$ is a TMS pair for all $g \in G$.
\item If $M$ is a compact normal subgroup of $G$, then $(KM/M,UM/M)$ is a TMS pair of $G/M$.
\item Let $V$ and $W$ be compact open subgroups of $G$ such that $V$ is contained in a conjugate of $UN$ for some compact normal subgroup $N$ of $G$, and $W \ge \bigcap_{g \in G}gUg\inv$.  Then the set of $g \in G$ such that $gKg\inv \le V$ but $gKg\inv \nleq W$ is compact and open, while the set of $g \in G$ such that $gKg\inv \le W$ does not have compact closure.
\end{enumerate}
\end{lemma}

\begin{proof}
Throughout the proof we refer to properties (a), (b) and (c) of the definition of a TMS pair.

We see by property (b) that $K$ is contained in a compact identity neighbourhood, so it has compact closure, and $K$ is infinite by property (c).  It is clear that the properties (a), (b) and (c) are invariant under taking conjugates and taking the closure, so any conjugate of $\ol{K}$ forms a TMS pair with $U$.  This proves (i).

Take a compact normal subgroup $M$ of $G$ and write $\widetilde{A}$ for $AM/M$, where $A$ is an element or subset of $G$.  It is clear that properties (b) and (c) are satisfied by $(\widetilde{K},\widetilde{U})$ in $\widetilde{G}$.  For property (a), take a compact normal subgroup $\widetilde{N} = N/M$ of $\widetilde{G}$ and a $\widetilde{G}$-conjugate $\widetilde{V} = \widetilde{h}\widetilde{U}\widetilde{h}\inv$ of $\widetilde{U}$, for some $h \in G$; note that $N$ is a compact normal subgroup of $G$ and $V := hUh\inv$ is a $G$-conjugate of $U$.  Given $g \in G$ such that $\widetilde{g}\widetilde{K}\widetilde{g}\inv$ is contained in $\widetilde{U}\widetilde{N}$ but not in $\widetilde{V}$, then $gKg\inv \le UN$ but $gKg\inv \nleq VN$; in particular, $gKg\inv \nleq V$.  By property (a), the set of elements $g \in G$ satisfying these conditions has compact closure, and hence the same is true of its image in $\widetilde{G}$.  Thus property (a) is satisfied by $(\widetilde{K},\widetilde{U})$, completing the proof of (ii).

For (iii), we first note that the set $H$ of $g \in G$ such that $gKg\inv \le V$ but $gKg\inv \nleq W$ is a union of right cosets of the compact open subgroup $V \cap W$; thus $H$ is clopen, and to show $H$ is compact it suffices to show it has compact closure.  By compactness we see that $W \ge \bigcap^n_{i=1}g_iUg\inv_i$, for some finite subset $\{g_1,\dots,g_n\}$ of $G$, and hence by properties (a) and (b) of a TMS pair, the set of $g \in G$ such that $gKg\inv \le W$ does not have compact closure.  To prove $H$ has compact closure, it is enough to show for $1 \le i \le n$ that the set of $g \in G$ such that $gKg\inv \le V$ but $gKg\inv \nleq g_iUg\inv_i$ has compact closure.  Thus we may assume $W$ is a conjugate of $U$.  We also have $V \le hUNh\inv$ for some $h \in G$, so that $H$ is contained in the set of $g \in G$ such that $gKg\inv \le hUN h\inv$ but $gKg\inv \nleq W$; thus we may assume $V = hUNh\inv$.  Finally, we see from property (a) that $h\inv H$ has compact closure (noting that $h\inv Wh$ is a conjugate of $U$), and hence $H$ has compact closure as required, completing the proof of (iii).
\end{proof}

Note that given Lemma~\ref{lem:tms_technicalities}(i), in defining a group with a TMS pair $(K,U)$, it makes little difference whether or not we require $K$ to be compact.

\begin{lemma}\label{lem:tms_separated}
Let $G$ be a compactly generated \tdlc group with a subgroup $K$ and let $X$ be a geometric $G$-set with some $G$-metric, such that $G$ acts faithfully.  Suppose $(K,U)$ is a TMS pair for some compact open subgroup $U$ such that $\bigcap_{g \in G}gUg\inv = \triv$.
\begin{enumerate}[label=(\roman*)]
\item $X$ is not quasi-isometric to a line.
\item For all $r \ge 0$ sufficiently large, writing $X^{K,r}$ for the set of points in $x \in X$ such that $K$ fixes pointwise the ball of radius $r$ around $x$, then $X^{K,r}$ is coarsely almost separated, and both $X^{K,r}$ and $X \setminus X^{K,r}$ are infinite.  If $X = G/V$ for some $V \le U$, we can take all $r \ge 0$.
\end{enumerate}
\end{lemma}

\begin{proof}
Since $K$ has the same fixed points as $\ol{K}$, we may assume that $K$ is closed, hence compact.
  
For (i), suppose $X$ is quasi-isometric to a line.  Then $G$ has a compact open normal subgroup $N$ such that $G/N$ is cyclic or dihedral.  But then $|K:N \cap K|$ is finite, a contradiction.

For (ii), we first consider the case that $X = G/V$ for some open $V \le U$, and the $G$-metric is given by a Cayley--Abels graph $\Gamma$ with vertex set $G/V$ and edges of $\Gamma$  given by $(gV,gsV)$ for $s \in S$, where $S$ is a compact symmetric generating set such that $S = VSV$.  Note that $(K,V)$ is a TMS pair of $G$ by Lemma~\ref{lem:tms_technicalities}(iii).

Let $X^K$ be the set of fixed points of $K$ on $X$ and let $V_1 = \bigcap_{s \in S}sVs\inv$; since $S$ is compact, $V_1$ is a finite intersection of conjugates of $V$.  Given $gV \in X^K$ with a neighbouring vertex $w \in X \setminus X^K$, then $w = gsV$ for some $s \in S$.  After conjugating, we find that $g\inv K g \le V$, but $g\inv K g \nleq V_1$.  Since $(K,V)$ is a TMS pair, this means $g$ is confined to a compact set, leaving only finitely many possibilities for $gV$ and hence for $w$.  We deduce that $X^K$ is an almost separated set of $X$.

Suppose for a contradiction that $X^K$ is finite, say $X^K = \{g_1V,\dots,g_nV\}$; and let $k \in K \setminus \triv$.  We can then find a compact open subgroup $W$ that is disjoint from $\{g\inv_1 k g_1, \dots, g\inv_n k g_n\}$; since $\bigcap_{g \in G}gVg\inv = \triv$, in fact we can take $W$ to be a finite intersection of $G$-conjugates of $V$, including $V$ itself.  It then follows that no $G$-conjugate of $W$ contains $K$, so $W$ contains no $G$-conjugate of $K$, which is incompatible with property (b) of the TMS pair $(K,W)$.  Similarly, suppose for a contradiction that $X \setminus X^K$ is finite, say $X = X^K \cup \{g_1V,\dots,g_nV\}$.  Then the group $K' = K \cap \bigcap^n_{i=1} g_i V g\inv_i$ is of finite index in $K$, hence nontrivial, but fixes every point in $X$.  To put this another way, we have
\[
N := \cgrp{gK'g\inv \mid g \in G} \le V \le U,
\]
which contradicts the hypothesis $\bigcap_{g \in G}gUg\inv=\triv$.  From these contradictions, we conclude that $X^K$ is infinite with infinite complement.  Since $X^K$ is almost separated, for any given $r \ge 0$, we see that $X^{K,r}$ is a cofinite subset of $X^K$; thus $X^{K,r}$ is almost separated and both $X^{K,r}$ and $X \setminus X^{K,r}$ are infinite.

In the general case, $G$ acts on $X$ with finitely many orbits and compact open stabilisers.  Choose some point $x \in X$; there is then some $r_0$ such that for all $r \ge r_0$ and $x \in X$, the pointwise fixator $G_{x,r}$ of the ball of radius $r$ around $x$ is contained in $U$.  From the fact that the set of fixed points of $K$ on $G/G_{x,r}$ is infinite with infinite complement, we see that $X^{K,r}$ is infinite with infinite complement.  Let $r_1$ be large enough that every ball of radius $r_1$ intersects all orbits of $G$ and let $r \ge r_0+r_1$.  Fix $n \ge 0$ and suppose we have $y \in X^{K,r}$ and $z \in X \setminus X^{K,r}$ with $d_X(y,z) \le n$.  Then there are some points $y'$ and $z'$ within distance $r_1$ of $y$ and $z$ respectively that belong to $Gx$; say $y' = g_yx$ and $z' = g_zx$.  Since $y'$ and $z'$ are within a bounded distance of each other we can take $g_z = g_yg$ where $g$ is confined to a finite set $F_1$, chosen independently of $(y,z)$.  There is then a finite subset $F_2$ of $G$ such that $\bigcap_{f \in F_2}fG_{x,r-r_1}f\inv$ fixes all points within distance $r+r_1$ of $hx$ for all $h \in F_1$.

We now have $K \le G_{y',r-r_1}$, so $g\inv_y Kg_y \le G_{x,r-r_1}$, but $g\inv_y Kg_y$ moves a point within distance $r+r_1$ of $g\inv_y z' = gx$, so $g\inv_y Kg_y \nleq fG_{x,r-r_1}f\inv$ for some $f \in F_2$.  Since $(K,G_{x,r-r_1})$ is a TMS pair and $F_2$ is finite, we see that $g\inv_y$ is confined to a compact set, independently of the pair $(y,z)$.  Thus $z'$ is confined to a bounded set, and hence also $z$ is confined to a bounded set, showing that $\delta_n X^{K,r}$ is bounded.  Thus $X^{K,r}$ is coarsely almost separated.
\end{proof}

We next consider the scale group case.

\begin{lemma}\label{lemma:tms_end}
Let $G$ be a scale group with a TMS subgroup $K$, acting geometrically on the locally finite leafless tree $T$ with fixed end $\xi$.   Then $K$ fixes pointwise a half-tree $T'$ such that $\xi \in \Ends T'$.  Consequently, the action of $G$ on $\Ends T$ is nondiscretely micro-supported.
\end{lemma}

\begin{proof}
We may quotient out by the kernel of the action on $T$, and so assume $G$ acts faithfully on $T$.  It is then easy to see that $G$ has no nontrivial compact normal subgroups.  For each $v \in VT$ we see that the conjugacy class of $G_v$ forms a base of neighbourhoods of the identity in $G$; arguing as in Lemma~\ref{lem:tms_technicalities}(iii), we see that $(K,G_v)$ is a TMS pair.

Let $T^K$ be the fixed subtree of $K$.  We have $G = P \rtimes \grp{l}$, where $l$ is a translation towards $\xi$ and the $P$-orbits on $VT$ are exactly the horospheres around $\xi$.  Since $K$ is compact, it fixes pointwise a ray representing $\xi$, and hence $\xi \in \Ends T^K$.  By Lemma~\ref{lem:tms_separated}(ii) and the fact that $T^K$ is a subtree, we see that $VT^K$ is an almost separated subset of $VT$ with infinite complement.  We can thus take a half-tree $T' = T_a$, with $a$ on the axis of $l$ and pointing towards $\xi$, such that $t(a) \in VT^K$ and such that no vertex of $T_a$ belongs to $\delta_1(VT^K)^c$, ensuring that $T_a \subseteq T^K$; in other words, the support of $K$ is confined to the half-tree $T_{\ol{a}}$.  In particular, in the action on $\Ends T$, the support of $K$ is confined to a compact subset of $\Ends T \setminus \{\xi\}$.  Since $G$ has extremely proximal action on $\Ends T \setminus \{\xi\}$, it follows that every nonempty open subset of $\Ends T \setminus \{\xi\}$, and hence also of $\Ends T$, contains the support of some $G$-conjugate of $K$.  Since $K$ is nondiscrete and $G$ acts faithfully on $\Ends T$, we deduce that the action is nondiscretely micro-supported.
\end{proof}

We now conclude the proof of the theorem.

\begin{proof}[Proof of Theorem~\ref{cpctend}]
Choose a compact open subgroup $U/M$ of the quotient group $G/M$ such that $(K/M,U/M)$ is a TMS pair for $G/M$, let $\Gamma$ be a Cayley--Abels graph for $G/M$ with vertex set $G/U$, and let $\Gamma^K$ be the graph of fixed points of $K$ acting on $\Gamma$.  Without loss of generality, $M = \bigcap_{g \in G}gUg\inv$.  Then by Lemma~\ref{lem:tms_separated}(ii), there is an almost separated subset of $V\Gamma$ that is infinite with infinite complement.  Moreover, by Lemma~\ref{lem:tms_separated}(i), $G/M$ does not act geometrically on a line.  We deduce that $G/M$ acts arc-geometrically on a leafless tree $T$ that is not a line, and we have an action as in case (ii) or (iii) of Proposition~\ref{prop:tree_dense}, with compact open arc stabilisers.  In particular, the kernel $N/M$ of this action is compact, and on the other hand any compact normal subgroup of $G$ has bounded orbits, so is contained in $N$.  Thus $N$ is the largest compact normal subgroup of $G$, and we note that $\Ends G \cong \Ends G/N$ as $G$-spaces.  From now on we can pass to the quotient $G/N$, so we assume $G$ has no nontrivial compact normal subgroups; in particular, $M = \triv$.  We also note via Lemma~\ref{lem:tms_technicalities}(iii) that $(K,V)$ is a TMS pair for every open subgroup $V$ contained in a conjugate of $U$.

If $G$ is a scale group, the action of $G$ on $\Ends G$ is nondiscretely micro-supported by Lemma~\ref{lemma:tms_end}, and hence by Theorem~\ref{thm:universal_micro}, $G$ is [A]-semisimple.  Thus we may assume for the rest of the proof that $G$ is not a scale group.  Since $G$ does not act arc-geometrically on a line, we then see by Proposition~\ref{prop:scale_group} that in every arc-geometric action of $G$ on a leafless tree, $G$ fixes no ends of the tree, so in fact $G$ is of general infinitely-ended type.

Starting from the Cayley--Abels graph $\Gamma$ for $G$, form the trees $T_i := T^{(t+i)}$ as in the Dicks--Dunwoody construction, starting from some $t \ge 1$ large enough that $T_i$ has more than one vertex for all $i \ge 0$.  By Lemma~\ref{lem:DD_subtree} it follows that the action of $G$ on $T_i$ is arc-geometric with no $G$-invariant subtree.  Equip $AT_i$ with a $G$-metric $d_i$ as in Lemma~\ref{lem:tree_geometric}.  Note that we are in case (iii) of Proposition~\ref{prop:tree_dense}, so $G$ has faithful extremely proximal action on the infinite perfect space $\ol{\Ends} T_i$.  

Consider the set $(AT_i)^{K,r}$ of $a \in AT_i$ such that $K$ fixes every $b \in AT_i$ with $d_i(a,b) \le r$.  By Lemma~\ref{lem:tms_separated}(ii), for $r \ge 0$ sufficiently large, the set $Y = (AT_i)^{K,r}$ is coarsely almost separated, infinite and has infinite complement in $AT_i$.  Considering the metric on $AT_i$, we see that the span of $Y$ contains a half-tree of $T_i$, with the result that the support of $K$ on $\ol{\Ends} T_i$ is not dense.  Since the action of $G$ on $\ol{\Ends} T_i$ is extremely proximal, it follows that every nonempty open subset of $\ol{\Ends} T_i$ contains the support of some $G$-conjugate of the nondiscrete subgroup $K$; thus the action of $G$ on $\ol{\Ends} T_i$ is nondiscretely micro-supported.  In particular, $G$ is [A]-semisimple by Theorem~\ref{thm:universal_micro}.  The assertions of (ii)(a) are now clear; (ii)(b) follows by Proposition~\ref{prop:DD_approximation}.
\end{proof}

\subsection{Sources of TMS subgroups}

In light of Theorem~\ref{cpctend}, it is interesting to consider sufficient conditions for a compactly generated \tdlc group $G$ to have a TMS subgroup.  We first note a partial converse to Theorem~\ref{cpctend}.

\begin{prop}\label{prop:smooth_tree}
Let $G$ be a compactly generated \tdlc group with faithful arc-geometric action on a leafless tree $T$.  Let $K = \rist_G(T_a)$ for some arc $a \in AT$, and suppose that $K$ fixes only finitely many arcs of $T_a$.  Then $(K,G_a)$ is a TMS pair.
\end{prop}

\begin{proof}
Equip $AT$ with a $G$-metric.  The hypothesis ensures that $(AT)^K$ differs from $AT_{\ol{a}}$ by only finitely many elements; thus $(AT)^K$ is almost separated.  For property (a) of a TMS subgroup, it is enough to show the following: for any two arcs $b_1$ and $b_2$, the set $L_{b_1,b_2}$ of $g \in G$ such that $gKg\inv \le G_{b_1}$ but $gKg\inv \nleq G_{b_2}$ is compact.  So fix $b_1,b_2 \in AT$ and let $g \in L_{b_1,b_2}$.  Then we see that $K$ fixes $g\inv b_1$ but not $g\inv b_2$.  Thus $g\inv b_2 \in \delta_d (AT)^K$ where $d = d(b_1,b_2)$ does not depend on $g$, so $g\inv b_2$ is confined to a finite set, and hence $g$ is confined to a compact set.  Thus property (a) for a TMS pair is satisfied.

We now consider the type of action $G$ has on $T$; note that by Lemma~\ref{lem:tree_geometric}, $G$ preserves no proper subtree.  The fact that $K$ acts nontrivially means that $T$ cannot be a line; it is possible that $G$ fixes at most one end $\xi$, but since $K$ only fixes finitely many arcs of $T_a$, we see that $\xi$ cannot be an end of $T_a$.  In either case, we see by Proposition~\ref{prop:tree_dense} that $\ol{\Ends} T_a$ is a compressible subspace of $\ol{\Ends} T$ under the action of $G$.  From there it is easy to see that the set of $g \in G$ such that $gKg\inv \le G_a$ does not have compact closure, so property (b) is satisfied.  On the other hand, one sees that any compact normal subgroup $N$ of $G$ acts trivially on the tree, from which it follows that $|K:K \cap N|$ is infinite.  Thus property (c) is satisfied and we conclude that $(K,G_a)$ is a TMS pair.
\end{proof}

More interesting is to obtain a TMS subgroup without imposing much large-scale structure on $G$, but more from the local structure of $G$.  Given the conclusions of Theorem~\ref{cpctend}, a necessary condition is that $G$ has a largest compact normal subgroup $N$ and $G/N$ is faithful micro-supported.  For convenience let us take $N=\triv$, so we assume $G$ itself is faithful micro-supported and has no nontrivial compact normal subgroups.  In particular, via Theorem~\ref{thm:universal_micro}, $G$ is [A]-semisimple and acts faithfully on $\lcent(G)$.  In this context we find a sufficient condition for $K = \rist_G(\alpha)$ to be a TMS subgroup for $\alpha \in \lcent(G)$, with three variants depending on whether or not $\alpha \in \ldlat(G)$ and whether or not $G$ is locally of finite quotient type.

\begin{prop}\label{prop:arborescent}
Let $G$ be an [A]-semisimple \tdlc group without nontrivial compact normal subgroups and let $\mc{Q} \subseteq \lcent(G) \setminus \{0,\infty\}$ be nonempty and $G$-invariant, generating a subalgebra $\mc{B} = \grp{\mc{Q}}$.  Fix $\alpha \in \mc{Q}$ and let $K = \rist_G(\alpha)$, with all rigid stabilisers defined with respect to the action on $\mc{B}$.  Suppose that $\alpha$ and $\mc{Q}$ satisfy the following conditions:
\begin{enumerate}[label=(\roman*)]
\item Given $\beta \in \mc{B}$ such that $\beta > 0$, there is $0 < \gamma \le \beta$ such that $\gamma \in \mc{Q}$.
\item There is $g \in G$ such that $g\alpha < \alpha$.
\item For all $\beta \in \mc{Q}$, its stabiliser $G_{\beta}$ is compact.
\item At least one of the following holds:
\begin{enumerate}[label=(\arabic*)]
\item For each natural number $n$, there are only finitely many $\beta \in \mc{Q}$ such that $\beta < \alpha$ and
\[
|K:K_\beta| \le n.
\]
\item We have $\mc{Q} \subseteq \ldlat(G)$, and for each $\alpha \in \mc{Q}$ and natural number $n$, there are only finitely many $\beta \in \mc{Q}$ such that $\beta < \alpha$ and
\[
|K:\rist_G(\beta)\CC_K(\rist_G(\beta))| \le n.
\]
\item We have $\mc{Q} \subseteq \ldlat(G)$ and $G$ is locally of finite quotient type.
\end{enumerate}
\end{enumerate}
Then $(K,U)$ is a TMS pair of $G$ for every compact open subgroup $U$ of $G$.
\end{prop}

\begin{proof}
From (iii), together with the fact that $G$ has no nontrivial compact normal subgroups, we see that $G$ acts faithfully on $\mc{B}$.  By \cite[Proposition~5.16]{CRWpart1}, for all $\alpha \in \mc{B}$ the group $\rist_G(\alpha)$ is the largest locally normal subgroup representing $\alpha$; in particular, given distinct elements $\alpha$ and $\beta$ of $\mc{B}$ then $\rist_G(\alpha)$ and $\rist_G(\beta)$ do not have any open subgroup in common.  Note also that $\N_G(\rist_G(\alpha)) = G_\alpha$ for all $\alpha \in \mc{B}$; in particular, $\N_G(\rist_G(\alpha))$ is compact for all $\alpha \in \mc{Q}$.

Fix compact open subgroups $U$ and $V$ of $G$; we claim that the set $H = H_{U,V}$ of $h \in G$ such that $hKh\inv \le U$ but $hKh\inv \nleq V$ has compact closure.  We may assume that $U \nleq V$.  Since $U$ is compact and acts faithfully on $\mc{B}$ with finite orbits, there is a base of neighbourhoods of the identity in $U$ consisting of pointwise stabilisers of finite $U$-invariant subalgebras of $\mc{B}$.  Since $U \cap V$ is a neighbourhood of the identity in $U$, we infer that there is a finite $U$-invariant subalgebra $\mc{B}'$ of $\mc{B}$ whose pointwise stabiliser $W$ in $U$ is contained in $V$; note that $W$ is also the pointwise stabiliser in $U$ of the set $\mc{P} = \{\beta_1,\dots,\beta_m\}$ of atoms of $\mc{B}'$.  By condition (i), for each $1 \le i \le m$ there is $0 < \gamma_i \le \beta_i$ such that $\gamma_i \in \mc{Q}$.  Given $1 \le i \le m$, write $H_i$ for the set of $g \in G$ such that $gKg\inv \le U$ and $\alpha > g\inv\gamma_i$.  Let $n \in \bN$ be such that $|U:L_i| \le n$ for $1 \le i \le m$, where under hypothesis (1) we set $L_i = U_{\gamma_i}$ and under hypothesis (2) we set
\[
L_i = \rist_U(\gamma_i) \CC_U(\rist_U(\gamma_i)).
\]
Now consider $g \in H_i$.

Under hypothesis (1), the stabiliser of $\gamma_i$ in $gKg\inv$ has index at most $n$, so the stabiliser of $g\inv \gamma_i$ in $K$ has index at most $n$.  We deduce that there are only finitely many possibilities for $g\inv \gamma_i$.

Under hypothesis (2), we have $U \ge gKg\inv > \rist_G(\gamma_i)$, so
\[
\rist_G(\gamma_i) = \rist_{gKg\inv}(\gamma_i) = \rist_U(\gamma_i).
\]
After conjugating by $g$ we see that $\rist_G(g\inv\gamma_i) = \rist_K(g\inv\gamma_i)$ and
\[
|K: \rist_G(g\inv \gamma_i) \CC_K(\rist_G(g\inv \gamma_i))| = |K:K \cap g\inv L_ig| \le n.
\]
Hypothesis (2) then leaves only finitely many possibilities for $g\inv\gamma_i$.

Under hypothesis (3), let $W = K \times \rist_U(\alpha^\perp)$.  Given $\beta \in \mc{Q}$ and a natural number $n'$ such that $\beta < \alpha$ and 
\[
|K:\rist_G(\beta)\CC_K(\rist_G(\beta))| \le n',
\]
then $\rist_G(\beta)$ is a direct factor of an open subgroup of $W$ of index at most $n'$.  By Lemma~\ref{lem:strong_fin_decomp} this leaves only finitely many possibilities for $\beta$, so in fact, hypothesis (2) holds.

In all cases, since $G_{\gamma_i}$ is compact, we deduce that $H_i$ is compact.

Given $h \in H$, then $hKh\inv$ acts nontrivially on $\mc{P}$, so $\beta_i$ and $k\beta_i$ are disjoint for some $k \in hKh\inv$ and $1 \le i \le m$, and hence $h\alpha > \beta_i \ge \gamma_i$.  We deduce that 
\[
H \subseteq \bigcup^m_{i=1}H_i;
\]
thus $H$ has compact closure as claimed.  In particular, $(K,U)$ satisfies condition (a) for a TMS pair.

By hypothesis, there is $g \in G$ such that $g\alpha < \alpha$.  It then follows that $g$ is not contained in any compact subgroup and $g^nKg^{-n} < K$ for all $n \ge 0$.  Taking $W$ to be a compact open subgroup of $G$ containing $K$, it follows that the set of $h \in G$ such that $hKh\inv \le W$ is unbounded; since $H_{W,U}$ has compact closure, the set of $h \in G$ such that $hKh\inv \le U$ is also unbounded.  Thus $(K,U)$ satisfies condition (b) for a TMS pair.  

Finally, by hypothesis $G$ has no finite locally normal subgroups, so $K$ is infinite.  Since $G$ has no nontrivial compact normal subgroups, condition (c) for a TMS subgroup is immediate.  Thus $(K,U)$ is a TMS pair of $G$.
\end{proof}

To justify some of the hypotheses in Proposition~\ref{prop:arborescent}, we recall the following from \cite{CRWpart2}.

\begin{lemma}[See {\cite[Theorem~6.19]{CRWpart2}}]\label{lem:skewering}
Let $G$ be an [A]-semisimple \tdlc group without nontrivial compact normal subgroups and let $\mc{A}$ be a subalgebra of $\ldlat(G)$ on which $G$ acts faithfully.  Then there exist nonzero elements $\alpha_1,\dots,\alpha_d$ of $\mc{A}$, where $d$ is the number of minimal nontrivial closed normal subgroups of $G$, such that for all $\beta \in \mc{A} \setminus \{0\}$, there is $g \in G$ and $1 \le i \le d$ such that $g\alpha_i < \beta$.  In particular, setting $\mc{Q} = \{g\alpha_i \mid g \in G, 1 \le i \le d\}$ and $\alpha = \alpha_1$, then $\alpha$ and $\mc{Q}$ satisfy conditions (a) and (b) of Proposition~\ref{prop:arborescent}.
\end{lemma}

The critical conditions in Proposition~\ref{prop:arborescent} are thus (c) and (d).  In particular, we can now prove Theorem~\ref{intro:cpctend:lqft}.

\begin{proof}[Proof of Theorem~\ref{intro:cpctend:lqft}]
Since $G$ has no nontrivial compact normal subgroups, we for every $\alpha \in \mc{Q}$ that $G_{\alpha}$ contains no nontrivial normal subgroup of $G$, so $G$ acts faithfully on $\mc{B} = \grp{\mc{Q}}$.  By Theorem~\ref{thm:universal_micro}, $G$ is [A]-semisimple.  By our present hypotheses, conditions (a), (c) and (d)(3) of Proposition~\ref{prop:arborescent} are satisfied, and by Lemma~\ref{lem:skewering} we can take $\alpha \in \mc{Q}$ such that condition (b) is satisfied.  Thus $K = \rist_G(\alpha)$ is a TMS subgroup of $G$.

By \cite[Corollary~1.4]{ReidEndo} it follows that given a compact open subgroup $U$ of $G$, then $U$ is not isomorphic to any of its proper open subgroups.  In particular, we see that if $G$ acts on a tree with more than two ends and with compact open stabilisers, then $G$ cannot fix any end of the tree.    Thus we are in case (ii) of Theorem~\ref{intro:cpctend}.
\end{proof}

Theorem~\ref{intro:cpctend:lqft} can be used to give a sufficient condition for a one-ended compactly generated \tdlc group $G$ to be locally indecomposable.  The only non-obvious aspect of applying Theorem~\ref{intro:cpctend:lqft} in this manner is in imposing conditions on locally normal subgroups (without direct reference to the decomposition lattice) ensuring that elements of $\ldlat(G) \setminus \{0,\infty\}$ have compact stabilisers.

\begin{corollary}\label{cor:cpctend:bis}
Let $G$ be a nontrivial compactly generated \tdlc group.  Suppose the following: $G$ is one-ended; $\QZ(G)=\triv$; $G$ has no nontrivial compact normal subgroups; $G$ is locally of finite quotient type; no open subgroup of $G$ has an infinite discrete quotient; and the centraliser of every nontrivial closed locally normal subgroup is compact.  Then $G$ is locally indecomposable.
\end{corollary}

\begin{proof}
Suppose that $G$ is locally decomposable, that is, there is some $\alpha \in \ldlat(G) \setminus \{0,\infty\}$.  Then $G_{\alpha}$ is an open subgroup, such that an open normal subgroup of $G_{\alpha}$ splits as a direct product $H = \rist_G(\alpha) \times \rist_G(\alpha^c)$.  Since the groups $\rist_G(\alpha)$ and $\rist_G(\alpha^c)$ are closed locally normal subgroups that centralise each other, they are both compact, and hence $H$ is compact.  Since $G_{\alpha}$ has no infinite discrete quotient, $G_{\alpha}/H$ is finite, so $G_{\alpha}$ is also compact.

Our goal is to invoke Theorem~\ref{intro:cpctend:lqft} with $\mc{Q} = \ldlat(G) \setminus \{0,\infty\}$.  Since $G_\alpha$ is compact for all $\alpha \in \mc{Q}$ and $G$ has no nontrivial compact normal subgroup, it follows that the $G$-action on the Boolean algebra $\mc{B} = \grp{\mc{Q}}$ is faithful. Therefore $G$ is   [A]-semisimple by Theorem~\ref{thm:universal_micro}, and it follows from Lemma~\ref{lem:skewering} that all hypotheses of Theorem~\ref{intro:cpctend:lqft} are satisfied. It follows that $G$ has more than one end, contradicting our hypothesis.  So we must in fact have $\ldlat(G) = \{0,\infty\}$.
\end{proof}

The authors do not know a general method for proving that, given a \tdlc group $G$, that the centraliser of every nontrivial closed locally normal subgroup of $G$ is compact.  However, such a restriction on centralisers can be proved for many complete geometric Kac--Moody groups, allowing the application of Corollary~\ref{cor:cpctend:bis}; this is done in another article, \cite{LocindecKM}.

\bibliographystyle{abbrv} 
\bibliography{Micro-supported-tree_arXiv_published_ver}

\end{document}